\documentclass[a4paper,reqno]{amsart}

\usepackage{amssymb}
\usepackage{latexsym}
\usepackage{amsmath}
\usepackage{amsthm}
\usepackage{bbm}
\usepackage{pgfplots}
\pgfplotsset{compat=1.15}
\usepackage{mathrsfs}
\usetikzlibrary{arrows}
\definecolor{ududff}{rgb}{0.30196078431372547,0.30196078431372547,1}
\definecolor{zzttqq}{rgb}{0.6,0.2,0}
\definecolor{xdxdff}{rgb}{0.49019607843137253,0.49019607843137253,1}
\definecolor{ududff}{rgb}{0.30196078431372547,0.30196078431372547,1}
\usepackage{hyperref}
\usepackage{xcolor}
\usepackage{amsmath}
\def\sa{{\mathfrak a}} 
\def\phi{\varphi}

\def\cA{{\mathcal A}}

   \def\cN{{\mathcal N}}

\DeclareMathOperator\sgn{{\text{\rm sgn}}} 
\DeclareMathOperator\dom{{\text{\rm dom\,}}} %Domain

\DeclareMathOperator{\spann}{span}
\DeclareMathOperator{\N}{\mathbb{N}}
\DeclareMathOperator{\R}{\mathbb{R}}
\DeclareMathOperator{\C}{\mathbb{C}}
\DeclareMathOperator{\diver}{div}

\renewcommand\setminus{\mathbin{\mathpalette\rsetminusaux\relax}}
\newcommand\rsetminusaux[2]{\mspace{-4mu}
  \raisebox{\rsmraise{#1}\depth}{\rotatebox[origin=c]{-20}{$#1\smallsetminus$}}
 \mspace{-4mu}
}
\newcommand\rsmraise[1]{%
  \ifx#1\displaystyle .8\else
    \ifx#1\textstyle .8\else
      \ifx#1\scriptstyle .6\else
        .45%
      \fi
    \fi
  \fi}

\makeatletter
\newcommand*{\rom}[1]{\expandafter\@slowromancap\romannumeral #1@}
\makeatother

\newtheorem{theorem}{Theorem}[section]
\newtheorem*{thm*}{Theorem}
\newtheorem{proposition}[theorem]{Proposition}

\newtheorem{lemma}[theorem]{Lemma}

\theoremstyle{definition}

\theoremstyle{definition}

\newtheorem{remark}[theorem]{Remark}

\numberwithin{equation}{section}

\title[A strict Neumann--Dirichlet eigenvalue inequality]{A strict Neumann--Dirichlet eigenvalue inequality on planar domains}

\author[J.~Rohleder]{Jonathan Rohleder}
\address{Matematiska institutionen \\ Stockholms universitet \\
106 91 Stockholm \\
Sweden}
\email{jonathan.rohleder@math.su.se}

\begin{document}

\begin{abstract}
On simply connected, bounded planar Lipschitz domains it is known from our earlier work that the ($k+2$)-nd eigenvalue of the Neumann Laplacian cannot exceed the $k$-th Dirichlet eigenvalue. In this note we sharpen this statement by showing that this inequality is always strict. This was previously only known under additional assumptions.
\end{abstract}

\maketitle

\section{Introduction}

On a bounded, connected Lipschitz domain $\Omega \subset \R^d$, $d \geq 2$, let
\begin{align*}
 0 = \mu_1 < \mu_2 \leq \mu_3 \leq \dots
\end{align*}
denote the eigenvalues of the Laplacian $- \Delta_{\rm N}$ with Neumann boundary conditions and
\begin{align*}
 0 < \lambda_1 < \lambda_2 \leq \lambda_3 \leq \dots
\end{align*}
those of the Laplacian $- \Delta_{\rm D}$ with Dirichlet boundary conditions, all eigenvalues counted with multiplicities. The aim of this note is to prove the following theorem.

\begin{theorem}\label{thm:main}
If $\Omega \subset \R^2$ is a bounded, simply connected Lipschitz domain, then 
\begin{align}\label{eq:mainIneq}
 \mu_{k + 2} < \lambda_k
\end{align}
holds for all $k \in \N$.
\end{theorem}

The inequality $\mu_{k + 2} \leq \lambda_k$ for simply connected planar domains had before been established by the author \cite{R25}, including strict inequality if $\lambda_k$ is simple or $\partial \Omega$ contains a straight line segment. This confirmed partially the conjecture that, in dimension $d$, the inequality $\mu_{k + d} \leq \lambda_k$ should hold for all $k$, see, e.g., \cite[Conjecture 3.2.42]{LMP} or \cite{BLP09}. The proof of Theorem \ref{thm:main} given in this note is based on a non-trivial refinement of the mechanism in \cite{R25}; cf.\ Remark \ref{rem} below.

For convex, sufficiently smooth planar domains, $\mu_{k + 2} < \lambda_k$ is due to Payne \cite{P55}; its $d$-dimensional convex counterpart was established by Levine and Weinberger \cite{LW86}. For arbitrary bounded Lipschitz (and somewhat more general) domains in $\R^d$, the most general known result is the Friedlander--Filonov inequality, $\mu_{k + 1} < \lambda_k$ for all $k$ \cite{F91,F04}. Very recently some improvements for large $k$ \cite{FG26p} or large $d$ \cite{F26p} were considered. 

Finally, let us point out that the case $k = 1$ in \eqref{eq:mainIneq} has seen considerable progress quite recently. While $\mu_3 < \lambda_1$ on simply connected planar domains was proven by the author \cite{R25}, recently Tang and Zhang \cite{TZ26p} proved it for arbitrary bounded planar Lipschitz domains using a variant of the author's method combined with a certain isoperimetric inequality for an auxiliary problem; however, their proof seems not to carry over to higher $k$. Moreover, very recently Filonov \cite{F26p2} established $\mu_3 < \lambda_1$ in dimension $d \geq 3$ by using known isoperimetric inequalities for Dirichlet and Neumann eigenvalues.

\section{Preliminaries}

In this section we provide some preliminaries and recall briefly the main tool from \cite{R25}, the construction of a self-adjoint vector Laplacian whose spectrum coincides with the union of the non-zero Neumann and Dirichlet Laplacian eigenvalues if $\Omega$ is simply connected. 

Let $\Omega \subset \R^2$ be a bounded Lipschitz domain and denote by $\nu$ the exterior unit normal field defined almost everywhere on $\partial \Omega$. Recall that if $u \in L^2 (\Omega)^2$ is such that $\diver u \in L^2 (\Omega)$, then the normal trace $u |_{\partial \Omega} \cdot \nu \in H^{- 1/2} (\partial \Omega)$ is well-defined and
\begin{align*}
 \int_\Omega u \cdot \nabla \varphi + \int_\Omega \overline{\varphi} \diver u = \left(u |_{\partial \Omega} \cdot \nu, \varphi |_{\partial \Omega} \right)_{\partial \Omega}
\end{align*}
holds for all $\varphi \in H^1 (\Omega)$, where $\varphi |_{\partial \Omega} \in H^{1/2} (\partial \Omega)$ is the trace of $\varphi$ and $(\cdot, \cdot)_{\partial \Omega}$ denotes the sesquilinear duality between $H^{1/2} (\partial \Omega)$ and its dual $H^{-1/2} (\partial \Omega)$; see \cite[Chapter XIX, \S 1, Theorem 2]{DL}. Occasionally we make use of the unit tangent field $\tau = (- \nu_2, \nu_1)$ on $\partial \Omega$.

In the following we will also make use of the notation
\begin{align*}
 \nabla^\perp \varphi := \binom{- \partial_2 \varphi}{\partial_1 \varphi} \quad \text{and} \quad \omega (u) := \partial_1 u_2 - \partial_2 u_1
\end{align*}
if $\varphi \in H^1 (\Omega)$ and $u = (u_1, u_2)^\top \in L^2 (\Omega)^2$ is such that $\partial_1 u_2 - \partial_2 u_1$, taken distributionally, belongs to $L^2 (\Omega)$. More generally, if $u = (u_1, u_2)^\top$ is a vector field, we write $u^\perp := (- u_2, u_1)^\top$ for its pointwise orthogonal field obtained by a rotation of angle $\pi/2$ in positive direction.

Consider the sesquilinear form
\begin{align*}
 \sa [u, v] & = \int_\Omega \left( \diver u \, \overline{\diver v} + \omega (u) \overline{\omega (v)} \right)
\end{align*}
defined on
\begin{align*}
 \dom \sa & = \left\{ u \in L^2 (\Omega)^2 : \diver u, \omega (u) \in L^2 (\Omega), \nu \cdot u |_{\partial \Omega} = 0 \right\}.
\end{align*}
Using the theory of closed non-negative sesquilinear forms and associated self-adjoint operators (see, e.g., \cite[Chapter~VI]{K95}) it was shown in \cite[Proposition 3.1]{R25} that there exists a non-negative, self-adjoint operator $A$ in $L^2 (\Omega)^2$ such that
\begin{align*}
 \sa [u, v] = (A u, v)
\end{align*}
holds for all $u \in \dom A \subset \dom \sa$ and $v \in \dom \sa$. Furthermore, by the abstract theory, the domain of the operator $A$ consists of those $u \in \dom \sa$ for which there exists a vector field $w \in L^2 (\Omega)^2$ such that 
\begin{align*}
 \sa [u, v] = (w, v)
\end{align*}
holds for all $v \in \dom \sa$; in this case we have $A u = w$. Note that the operator $A$ has an interpretation as Hodge Laplacian with absolute boundary conditions; see \cite[Remark 3.9]{R25} and \cite{FGM25,HMZ26,M26}.

The following spectral description is a consequence of \cite[Proposition 3.2]{R25}.

\begin{proposition}\label{prop:spectrumA}
Assume that $\Omega \subset \R^2$ is a bounded, simply connected Lipschitz domain. Then the operator $A$ has purely discrete spectrum and for each $\lambda \in \R$ the identity
\begin{align}\label{eq:EVFormula}
 \ker (A - \lambda) = \nabla \ker (- \Delta_{\rm N} - \lambda) \oplus \nabla^\perp \ker (- \Delta_{\rm D} - \lambda)
\end{align}
holds.
\end{proposition}

In particular, the spectrum of $A$ consists of all positive eigenvalues of $- \Delta_{\rm N}$ and all eigenvalues of $- \Delta_{\rm D}$, appearing according to their multiplicities, if $\Omega$ is simply connected. Without the latter assumption, $A$ has in addition a non-trivial kernel, consisting of certain harmonic fields, whose dimension equals the first Betti number of $\Omega$; cf.\ \cite[Remark 3.7]{R25}.

Let us denote by $\eta_1 \leq \eta_2 \leq \dots$ the eigenvalues of $A$, again counted with multiplicities. Again by classical theory, we have
\begin{align}\label{eq:minMax}
 \eta_k = \min_{\substack{U \subset \dom \sa \\ \dim U = k}} \max_{\substack{u \in U \\ u \neq 0}} \frac{\int_\Omega \left( |\diver u|^2 + |\omega (u)|^2 \right)}{\int_\Omega |u|^2}
\end{align}
for each $k \in \N$.

\section{Proof of Theorem \ref{thm:main}}

In this section we prove Theorem \ref{thm:main}. We first collect a few useful lemmas. We assume throughout that $\Omega \subset \R^2$ is a bounded, simply connected Lipschitz domain.

\begin{lemma}\label{lem:goodProperties}
The following assertions hold true.
\begin{enumerate}
 \item If $\lambda > 0$ and $u \in \ker (A - \lambda)$, then 
\begin{align*}
 \omega (u) \in \ker (- \Delta_{\rm D} - \lambda) \subset H_0^1 (\Omega) \quad \text{and} \quad \diver u \in \ker (- \Delta_{\rm N} - \lambda) \subset H^1 (\Omega).
\end{align*}
 \item If $u, u^\perp \in \ker (A - \lambda)$, then $u = 0$ identically.
\end{enumerate}
\end{lemma}

\begin{proof}
By Proposition \ref{prop:spectrumA}, each $u \in \ker (A - \lambda)$ can be written
\begin{align}\label{eq:decomp}
 u = \nabla \psi + \nabla^\perp \varphi
\end{align}
for some $\psi \in \ker (- \Delta_{\rm N} - \lambda)$ and some $\varphi \in \ker (- \Delta_{\rm D} - \lambda)$. Hence,
\begin{align*}
 \omega (u) = \Delta \varphi = - \lambda \varphi \in \ker (- \Delta_{\rm D} - \lambda) \subset H_0^1 (\Omega)
\end{align*}
and
\begin{align}\label{eq:core}
 \diver u = \Delta \psi = - \lambda \psi \in \ker (- \Delta_{\rm N} - \lambda) \subset H^1 (\Omega);
\end{align}
this proves (i). 

If $u, u^\perp \in \ker (A - \lambda)$ and $u$ is written as in \eqref{eq:decomp}, then, by (i),
\begin{align*}
 H_0^1 (\Omega) \ni \omega (u^\perp) = \diver u = - \lambda \psi,
\end{align*}
so that $\psi \in \ker (- \Delta_{\rm N} - \lambda) \cap H_0^1 (\Omega)$; unique continuation (e.g.\ \cite[Lemma 4.4]{R25}) gives $\psi = 0$ constantly. However, in this case, $u = \nabla^\perp \varphi$ and, by (i),
\begin{align*}
 \lambda \varphi = - \Delta \varphi = \diver u^\perp \in \ker (- \Delta_{\rm N} - \lambda),
\end{align*}
that is, $\varphi \in H_0^1 (\Omega) \cap \ker (- \Delta_{\rm N} - \lambda)$. Another application of unique continuation yields $\varphi = 0$ identically and, thus, $u = 0$.
\end{proof}

In the following, for $\psi \in H^1 (\Omega)$ such that $\Delta \psi \in L^2 (\Omega)$ we denote by $\partial_\nu \psi |_{\partial \Omega} = \nabla \psi |_{\partial \Omega} \cdot \nu$ its derivative with respect to the exterior unit normal field on $\partial \Omega$.

\begin{lemma}\label{lem:curlIdentity}
If $u_1, u_2 \in H_0^1 (\Omega) \cap H^2 (\Omega)$ and $u = (u_1, u_2)^\top$, then 
\begin{align*}
 \omega (u) = \nu_1 \partial_\nu u_2 - \nu_2 \partial_\nu u_1
\end{align*}
holds almost everywhere on $\partial \Omega$.
\end{lemma}

\begin{proof}
Since $u_j \in H_0^1 (\Omega) \cap H^2 (\Omega)$, we have $\partial_i u_j \in H^1 (\Omega)$, $i = 1, 2$, and $\nabla u_j \cdot \tau = 0$ on $\partial \Omega$, where the components of $\nabla u_j$ on $\partial \Omega$ are traces of functions in $H^1 (\Omega)$ and, hence, belong to $L^2 (\partial \Omega)$. Consequently, we get
\begin{align*}
 \nabla u_j = (\nabla u_j \cdot \tau) \tau + (\nabla u_j \cdot \nu) \nu = (\partial_\nu u_j) \nu
\end{align*}
almost everywhere on $\partial \Omega$. In particular, $\partial_i u_j = \nu_i \partial_\nu u_j$ for $i = 1, 2$, and the claim follows from the definition of $\omega (u)$.
\end{proof}

\begin{lemma}\label{lem:curlZero}
If $u_1, u_2 \in \ker (- \Delta_{\rm D} - \lambda)$ and $\omega (u) \in H_0^1 (\Omega)$, where $u = (u_1, u_2)^\top$, then $u \in \ker (A - \lambda)$. 
\end{lemma}

\begin{proof}
As $u_1, u_2 \in H_0^1 (\Omega)$, $u$ belongs to $\dom \sa$. Furthermore, $\diver u \in H^1 (\Omega)$ due to
\begin{align*}
 \nabla \diver u = \Delta u - \nabla^\perp \omega (u) \in L^2 (\Omega)^2.
\end{align*}
Hence for each $v \in \dom \sa$ we have
\begin{align*}
 \sa [u, v] & = \int_\Omega \left( \diver u \, \overline{\diver v} + \omega (u) \overline{\omega (v)} \right) = - \int_\Omega \Delta u \cdot v = \lambda \int_\Omega u \cdot v,
\end{align*}
where the boundary terms vanish due to $v |_{\partial \Omega} \cdot \nu = 0$ and $\omega (u) \in H_0^1 (\Omega)$. This implies the assertion.
\end{proof}

The following proposition will be the main ingredient in the proof of the strict inequality in Theorem \ref{thm:main}.

\begin{proposition}\label{prop:mainStep}
Assume that $\lambda$ is an eigenvalue of $- \Delta_{\rm D}$. Moreover, assume that $W \subset \ker (- \Delta_{\rm D} - \lambda)$ is a subspace with $p := \dim W \geq 1$. Then
\begin{align*}
 \dim \left( W^2 \cap \ker (A - \lambda) \right) \leq p - 1
\end{align*}
holds.
\end{proposition}

\begin{proof}
Let us abbreviate
\begin{align*}
 D_W := W^2 \cap \ker (A - \lambda),
\end{align*}
so that our goal is to prove $\dim D_W \leq p - 1$. Assume for a contradiction that $\dim D_W \geq p$.

{\bf Step 1.} In this step we derive $W \subset H^2 (\Omega)$ from the contradiction assumption. Let $\{\varphi_1, \dots, \varphi_p\}$ be an orthonormal basis of $W$. By the assumption $\dim D_W \geq p$, $D_W$ contains $p$ linearly independent elements 
\begin{align*}
 u^j = \sum_{k = 1}^p A_k^j \varphi_k, \quad j = 1, \dots, p,
\end{align*}
with $A_k^j \in \C^2$. Define the vectors
\begin{align*}
 \cA^j = \begin{pmatrix} A_1^j \\ \vdots \\ A_p^j \end{pmatrix} \quad \text{and} \quad \cA^j_\perp = \begin{pmatrix} (A_1^j)^\perp \\ \vdots \\ (A_p^j)^\perp \end{pmatrix}, \quad j = 1, \dots, p
\end{align*}
in $\C^{2 p}$ and the matrix
\begin{align*}
 \cA = \begin{pmatrix} | & & | & | & & | \\ \cA^1 & \cdots & \cA^p & \cA^1_\perp & \cdots & \cA^p_\perp \\ | & & | & | & & |  \end{pmatrix} \in \C^{2p \times 2p}.
\end{align*}
The matrix $\cA$ is invertible: in fact, assume $\cA C = 0$ for some vector $C = (c_1, \dots, c_{2p})^\top \in \C^{2p}$. Then
\begin{align*}
 \sum_{j = 1}^p c_j \begin{pmatrix} A_1^j \varphi_1 \\ \vdots \\ A_p^j \varphi_p \end{pmatrix} + \sum_{j = 1}^p c_{j + p} \begin{pmatrix} (A_1^j)^\perp \varphi_1 \\ \vdots \\ (A_p^j)^\perp \varphi_p \end{pmatrix} = 0.
\end{align*}
This can be read as $p$ equations in $\C^2$; summing up these equations yields 
\begin{align*}
 v + w^\perp = 0, \quad \text{where} \quad v = \sum_{j = 1}^p c_j u^j, \quad w = \sum_{j = 1}^p c_{j+p} u^j.
\end{align*}
In particular, $v, w \in \ker (A - \lambda)$ and $v = - w^\perp$. From this and Lemma \ref{lem:goodProperties} (ii) we obtain $v = 0 = w$ identically. Moreover, as the $u^j$ are linearly independent, $c_1 = \dots = c_{2p} = 0$ follows. Hence $\cA$ is invertible.

Next, note that, by Lemma \ref{lem:goodProperties} (i),
\begin{align}\label{eq:boundary0}
 H_0^1 (\Omega) \ni - \omega (u^j) = \sum_{k = 1}^p (A_k^j)^\perp \cdot \nabla \varphi_k, \quad j = 1, \dots, p,
\end{align}
and
\begin{align*}
 H^1 (\Omega) \ni \diver (u^j) = \sum_{k = 1}^p A_k^j \cdot \nabla \varphi_k, \quad j = 1, \dots, p.
\end{align*}
In particular,
\begin{align*}
 (\partial_1 \varphi_1, \partial_2 \varphi_1, \partial_1 \varphi_2, \partial_2 \varphi_2, \dots, \partial_1 \varphi_p, \partial_2 \varphi_p) \cA \in H^1 (\Omega)^{2p},
\end{align*}
and multiplication by $\cA^{-1}$ from the right yields $\varphi_k \in H^2 (\Omega)$, $k = 1, \dots, p$. Thus $W \subset H^2 (\Omega)$.

{\bf Step 2.} Define 
\begin{align*}
 V_W := \left\{ \partial_\nu \varphi |_{\partial \Omega} : \varphi \in W \right\};
\end{align*}
note that unique continuation implies that the normal derivative, restricted to $\ker (- \Delta_{\rm D} - \lambda)$, is injective and, hence, $\dim V_W = \dim W = p$. Furthermore, let
\begin{align*}
 H := \left\{ h \in L^2 (\partial \Omega) : h \nu_1, h \nu_2 \in V_W \right\}.
\end{align*}
In this step we prove that 
\begin{align}\label{eq:dimH}
 \dim H \geq \dim D_W.
\end{align}
For this, define the mapping 
\begin{align*}
 \cN : D_W \to H, \quad \cN u = \nu_1 \partial_\nu u_1 |_{\partial \Omega} + \nu_2 \partial_\nu u_2 |_{\partial \Omega}.
\end{align*}
It is well-defined since for $u = (u_1, u_2)^\top \in D_W \subset W^2$, on $\partial \Omega$ we have
\begin{align*}
 \nu_1 \cN u & = \nu_1^2 \partial_\nu u_1 + \nu_1 \nu_2 \partial_\nu u_2 = \partial_\nu u_1 + \nu_2 \left( \nu_1 \partial_\nu u_2 - \nu_2 \partial_\nu u_1 \right) = \partial_\nu u_1 + \nu_2 \omega (u) \\ 
 & = \partial_\nu u_1,
\end{align*}
where we have used Lemma \ref{lem:goodProperties} (i) and Lemma \ref{lem:curlIdentity}; analogously, $\nu_2 \cN u = \partial_\nu u_2$. The mapping $\cN$ is injective: in fact, if $\cN u = 0$, then the vector field $\binom{\partial_\nu u_1}{\partial_\nu u_2}$ is orthogonal to the normal vector field on $\partial \Omega$. On the other hand, by Lemma \ref{lem:goodProperties} (i), $u \in D_W \subset \ker (A - \lambda)$ implies, on $\partial \Omega$,
\begin{align*}
 0 = \omega (u) = \nu_1 \partial_\nu u_2 - \nu_2 \partial_\nu u_1,
\end{align*}
so that $\binom{\partial_\nu u_1}{\partial_\nu u_2}$ is even orthogonal to the tangential vector field $\tau = \nu^\perp$ and, hence, vanishes on $\partial \Omega$. Since $u_1, u_2 \in \ker (- \Delta_{\rm D} - \lambda)$, unique continuation implies $u = 0$ in $\Omega$. Thus $\cN$ is injective and $\dim H \geq \dim D_W$ follows. In fact, it is not difficult to see that $\cN$ is actually bijective and the latter dimensions coincide; however, as this is not used in this proof, we omit the details.

% In order to show that $\cN$ is surjective, we explicitly construct its inverse. In fact, define the mapping
% \begin{align*}
%  \gamma_\nu : W \to V_W, \quad \gamma_\nu w = \partial_\nu w |_{\partial \Omega}.
% \end{align*}
% This mapping is bijective. Hence, for $h \in H$, we have 
% \begin{align*}
%  u := \binom{\gamma_\nu^{-1} (h \nu_1)}{\gamma_\nu^{-1} (h \nu_2)} \in W^2
% \end{align*}
% and $(\partial_\nu u_1, \partial_\nu u_2)^\top = h \nu$ on $\partial \Omega$; in particular, by Lemma \ref{lem:curlIdentity},
% \begin{align*}
%  \omega(u) = \nu_1 \partial_\nu u_2 - \nu_2 \partial_\nu u_1 = \nu_1 h \nu_2 - \nu_2 h \nu_1 = 0
% \end{align*}
% on $\partial \Omega$; together with Step~1, $\omega (u) \in H_0^1 (\Omega)$ follows. Thus, Lemma \ref{lem:curlZero} implies $u \in D_W$. Finally, $\cN u = h \nu_1^2 + h \nu_2^2 = h$. Hence, $\cN$ is bijective, and \eqref{eq:dimH} is proven.

{\bf Step 3.} We finally derive a contradiction. Consider the mappings $P_1, P_2 : H \to V_W$ given by
\begin{align*}
 P_j h = h \nu_j, \quad h \in H, \quad j = 1, 2.
\end{align*}
Each of these mappings is injective. In fact, as $\Omega$ is a Lipschitz domain, there exist non-empty open sets $\Sigma_1, \Sigma_2 \subset \partial \Omega$ such that
\begin{align*}
 \nu_j \neq 0~\text{almost everywhere on}~\Sigma_j, \quad j = 1, 2;
\end{align*}
see Lemma \ref{lem:lipschitzNormal}. Now assume $P_1 h = 0$ for some $h \in H$. Then, $h = 0$ almost everywhere on $\{ \nu_1 \neq 0\}$ and, in particular, $h = 0$ almost everywhere on $\Sigma_1$. Furthermore, as $h \nu_2 \in V_W$ by definition of $H$, there exists some $\varphi \in W$ such that $h \nu_2 = \partial_\nu \varphi$ on $\partial \Omega$. As the left-hand side vanishes on $\Sigma_1$, unique continuation gives $\varphi = 0$ identically and, hence, $h = 0$ almost everywhere on $\{\nu_2 \neq 0\}$. As $|\nu| = 1$ almost everywhere on $\partial \Omega$, we have
\begin{align*}
 |\{\nu_1 \neq 0\} \cup \{\nu_2 \neq 0\}| = |\partial \Omega|
\end{align*}
($|\cdot|$ denoting the 1-dimensional Lebesgue measure on $\partial \Omega$). Hence, $h = 0$ almost everywhere on $\partial \Omega$. Injectivity of $P_2$ is analogous. 

Now, by Step 2 and the initial assumption of this proof, $\dim H \geq \dim D_W \geq p$, and since $P_j : H \to V_W$ is injective for $j = 1, 2$, and $\dim V_W = p$, it follows that $P_1$ and $P_2$ are isomorphisms. Consequently, 
\begin{align*}
 B := P_2 \circ P_1^{-1} : V_W \to V_W 
\end{align*}
is a well-defined isomorphism. By surjectivity of $P_1$, each $f \in V_W$ can be written $f = h \nu_1$ for some $h \in H$, and 
\begin{align}\label{eq:B}
 B f = \big(P_2 \circ P_1^{-1} \big) (h \nu_1) = P_2 h = h \nu_2.
\end{align}
As isomorphism between finite-dimensional complex vector spaces, $B$ has an eigenvalue $z \in \C$, $z \neq 0$. Assume that $f = h \nu_1 \in V_W \setminus \{0\}$ is a corresponding eigenfunction. From \eqref{eq:B} we obtain
\begin{align*}
 z h \nu_1 = z f = B f = h \nu_2,
\end{align*}
or, equivalently, $(z \nu_1 - \nu_2) h = 0$ on $\partial \Omega$. The rest of the argument is similar to the injectivity of $P_1$: Namely, we have $h = 0$ almost everywhere on $\{z \nu_1 - \nu_2 \neq 0\}$, and this set contains a non-empty, open set $\Sigma$. Moreover, $\nu_2 h = \partial_\nu \varphi$ for some $\varphi \in W$, and the left-hand side vanishes on $\Sigma$, forcing $\varphi = 0$ by unique continuation. Hence, $h = 0$ almost everywhere on $\{ \nu_2 \neq 0\}$, and since at almost every point of $\partial \Omega$ either $z \nu_1 - \nu_2 \neq 0$ or $\nu_2 \neq 0$ (note that $z \neq 0$), $h = 0$ almost everywhere follows, a contradiction. The proof is complete.
\end{proof}

Now we are ready to prove the main result of this note.

\begin{proof}[Proof of Theorem \ref{thm:main}]
For convenience of notation we assume that $k \geq 2$; the case $k = 1$ was already treated in \cite{R25} using that $\lambda_1$ is simple. Let us choose an orthonormal basis $\varphi_1, \varphi_2, \dots$ of $L^2 (\Omega)$ such that $- \Delta_{\rm D} \varphi_j = \lambda_j \varphi_j$ holds for each $j$. Now, fix some $k$ and let 
\begin{align*}
 W := \spann \left\{ \varphi_j : j \leq k, \lambda_j = \lambda_k \right\}
\end{align*}
and $p := \dim W$; that is, $k$ is the $p$-th index whose associated eigenvalue equals $\lambda_k$. In particular, 
\begin{align*}
 \lambda_{k - p} < \lambda_{k - p + 1} = \dots = \lambda_k.
\end{align*}
(Note that $k - p$ is positive due to $k \geq 2$ and simplicity of $\lambda_1$.) Now let 
\begin{align}\label{eq:uvw}
 u = \sum_{j = 1}^{k - p} \binom{\alpha_j}{\beta_j} \varphi_j, \quad w = \sum_{j = k - p + 1}^k \binom{\alpha_j}{\beta_j} \varphi_j, \quad \text{and} \quad v \in \ker (A - \lambda_k)
\end{align}
for arbitrary $\alpha_j, \beta_j \in \C$; note that $w \in W^2$. The same computation as in the proof of \cite[Theorem 4.1]{R25} gives
\begin{align}\label{eq:computation}
 \sa [u + w + v] 
%  \leq \lambda_k \|u + w + v\|^2 - (\lambda_k - \lambda_{k - p}) \|u\|^2 
 \leq \lambda_k \|u + w + v\|^2.
\end{align}
It follows from Proposition \ref{prop:mainStep} that the dimension of the space spanned by arbitrary vectors of the form $u$, $v$, and $w$ of the form \eqref{eq:uvw} equals at least 
\begin{align*}
 2 (k - p) + \dim W^2 + \dim \ker (A - \lambda_k) \, - & \dim \big( W^2 \cap \ker (A - \lambda_k) \big) \\
  & \,\, \geq 2 k + \dim \ker (A - \lambda_k) - (p - 1),
\end{align*}
where we have used that $u = v + w$ for any such $u, v, w$ would imply $u = 0 = v + w$ since the right-hand side satisfies $- \Delta (v + w) = \lambda_k (v + w)$ componentwise, while the left-hand side does not, unless $u = 0$. By \eqref{eq:computation} and the min-max principle, the eigenvalues of $A$ satisfy
\begin{align*}
 \eta_{2 k + \dim \ker (A - \lambda_k) - (p - 1)} \leq \lambda_k.
\end{align*}
As at most $\dim \ker (A - \lambda_k)$ many of the $\eta_j$ equal $\lambda_k$, it follows
\begin{align*}
 \eta_{2 k - p + 1} < \lambda_k.
\end{align*}
Since $- \Delta_{\rm D}$ has precisely $k - p$ eigenvalues in $(0, \lambda_k)$, it follows that $- \Delta_{\rm N}$ has at least 
\begin{align*}
 2 k - p + 1 - (k - p) = k + 1
\end{align*}
eigenvalues in $(0, \lambda_k)$. Taking into account also its zero eigenvalue, we arrive at $\mu_{k + 2} < \lambda_k$. As $k$ was arbitrary, this is the assertion of the theorem.
\end{proof}

\begin{remark}\label{rem}
It can be seen from the proof of Theorem \ref{thm:main} that Proposition \ref{prop:mainStep} can be avoided whenever it is possible to prove that
\begin{align*}
 \spann \left\{ \binom{\varphi_j}{0}, \binom{0}{\varphi_j} : j \leq k \right\} \cap \ker (A - \lambda_k) = \{0\}.
\end{align*}
However, the latter may be false in certain cases. For instance, if $\Omega$ is the unit disk, then $\lambda_2 = \lambda_3 = j_{1,1}^2 =: \lambda$, where $j_{1,1}$ is the first positive zero of the Bessel function $J_1$. The corresponding eigenspace is spanned, in polar coordinates, by the functions $\varphi_2 (r, \theta) = J_1 (j_{1,1} r) \cos \theta$ and $\varphi_3 (r, \theta) = J_1 (j_{1,1} r) \sin \theta$. On the other hand, the function
\begin{align*}
 \psi (r) = J_0 (j_{1,1} r)
\end{align*}
belongs to $\ker (- \Delta_{\rm N} - \lambda)$ (since $J_0' = - J_1$) and $(\nabla \psi) (r, \theta) = - j_{1,1} J_1 (j_{1,1} r) \binom{\cos \theta}{\sin \theta}$. Thus
\begin{align*}
 \ker (A - \lambda) \ni \nabla \psi \in \spann \left\{ \binom{\varphi_j}{0}, \binom{0}{\varphi_j} : j \leq 3 \right\}.
\end{align*}
If we set $W = \ker (- \Delta_{\rm D} - \lambda)$ and $p = \dim W = 2$, then we obtain from the above and Proposition \ref{prop:mainStep} that $\dim D_W = 1 = p - 1$ in this case. In particular, the assertion of Proposition \ref{prop:mainStep} is sharp.
\end{remark}

\begin{remark}
Proposition \ref{prop:mainStep} generalizes the argument given in \cite{R25} to show that $\mu_{k + 2} < \lambda_k$ holds whenever $\lambda_k$ is simple in the spectrum of the Dirichlet Laplacian, and $\Omega$ is simply connected.
\end{remark}

\appendix

\section{An observation on Lipschitz domains}

The following property of Lipschitz domains is quite intuitive; for the sake of completeness we provide a proof.

\begin{lemma}\label{lem:lipschitzNormal}
Assume that $\Omega \subset \R^2$ is a bounded Lipschitz domain and let $S \in \R^2 \setminus \{0\}$. Then there exists a non-empty, open set $\Sigma \subset \partial \Omega$ such that $\nu \notin \spann \{S\}$ almost everywhere on $\Sigma$.
\end{lemma}

\begin{proof}
Possibly after rotating $\Omega$ suitably, it suffices to prove the existence of a non-empty, open set $\Sigma \subset \partial \Omega$ such that $\nu_1 \neq 0$ (i.e.\ $\nu \neq \pm \textbf{e}_2$) almost everywhere on $\Sigma$. Assume, for a contradiction, that each non-empty, open subset of $\partial \Omega$ contains a set of positive measure at every point of which the normal vector is $\pm \textbf{e}_2$.

{\bf Claim:} For each point $p$ of $\partial \Omega$ there exists an open neighborhood $\Sigma_p \subset \partial \Omega$ of $p$ and a number $\sigma_p \in \{\pm 1\}$ such that
\begin{align*}
 q \in \Sigma_p \,\, \text{and}\,\, \nu (q) = \pm \textbf{e}_2 \quad \Longrightarrow \quad \nu (q) = \sigma_p \textbf{e}_2.
\end{align*}
Namely, for each $p \in \partial \Omega$, by definition of a Lipschitz domain, there exists an open neighborhood $U_p \subset \R^2$ of $p$ such that $\Sigma_p := U_p \cap \partial \Omega$ is, up to rotation, the graph of a Lipschitz continuous function and $U_p \setminus \Sigma_p$ has two connected components, of which one is contained in $\Omega$ and the other has empty intersection with $\Omega$. In particular, if $U_p$ is chosen sufficiently small, there exists a unit vector $d_p$ such that 
\begin{align}\label{eq:nuDirection}
 \nu (q) \cdot d_p < 0 \quad \text{for almost all}~q \in \Sigma_p. 
\end{align}
By the contradiction assumption, $d_p \neq \pm \textbf{e}_1$. Thus, $\sigma_p := - \sgn (d_p \cdot \textbf{e}_2) \neq 0$. In particular, if $q \in \Sigma_p$ is such that $\nu (q) = \pm \textbf{e}_2$, then $\nu (q) = \sigma_p \textbf{e}_2$, which proves the claim.

Let now $p^\pm$ be two points of $\overline{\Omega}$ with maximal, respectively minimal, second coordinate. These exist as $\overline{\Omega}$ is compact. Furthermore, $p^\pm \in \partial \Omega$ and $p^\pm \in \Gamma_{\rm ext}$, where $\Gamma_{\rm ext}$ denotes the connected component of $\partial \Omega$ whose enclosed set contains $\Omega$. By the choice of $p^\pm$, we have $\sigma_{p^\pm} = \pm 1$. On the other hand, $\Gamma_{\rm ext}$ is a connected, compact curve and $p \mapsto \sigma_p$ is locally constant along $\Gamma_{\rm ext}$ by the above claim and the contradiction assumption, hence constant. This is a contradiction.
\end{proof}

\subsection*{AI disclosure}

The main mathematical ideas in this note were developed through discussions with AI tools. The text itself is written entirely by the author.

\subsection*{Acknowledgments}
The  author gratefully acknowledges financial support by the grant no. 2022-03342 of the Swedish Research Council (VR).

\end{document}